\documentstyle[11pt]{article}

\input amssym.def 
\input amssym

\newcommand{\be}{\begin{eqnarray*}}
\newcommand{\ee}{\end{eqnarray*}}

\newcommand{\al}{\alpha}

\newcommand{\lon}{\longrightarrow}
\newcommand{\complex}{{\Bbb C}}

\newcommand{\half}{\frac{1}{2}}
\newcommand{\third}{\frac{1}{3}}

\newtheorem{thm}{Theorem}[section]

\newtheorem{cor}[thm]{Corollary}
\newtheorem{pro}[thm]{Proposition}

\newtheorem{defi}[thm]{Definition}

\newcommand{\frakg}{{\frak g}}
\newcommand{\frakh}{{\frak h}}

\newcommand{\calx}{{\cal X}}
\newcommand{\call}{{\cal L}}
\newcommand{\gm}{\Gamma }

\newcommand{\Alt}{\mbox{Alt}}
\newcommand{\smalcirc}{\mbox{\tiny{$\circ $}}}
\newcommand{\h}{{\frak h}}
\newcommand{\g}{{\frak g}}
\newcommand{\xa}{\xi_{\al}}
\newcommand{\ax}{\xi_{-\al}}
\newcommand{\calp}{{\cal P}}

\begin{document}

\title{{\bf Spin Calogero-Moser systems associated
with simple Lie algebras}\\ (In memory of J\"urgen 
 Moser)}

\author{Luen-Chau Li\\
{\sf email: luenli@math.psu.edu }\\
AND\\
Ping Xu
         \thanks{ Research partially supported by NSF
       grant DMS97-04391. }\\
{\sf email: ping@math.psu.edu }\\
        Department of Mathematics\\
         Pennsylvania State University \\
         University Park, PA 16802, USA}

\date{}
\maketitle
\begin{abstract}
We introduce spin Calogero-Moser systems associated with root systems
of simple Lie algebras and give the associated Lax representations
 (with spectral parameter) and fundamental Poisson
 bracket relations. Our analysis is based
on a group-theoretic framework similar in spirit to the standard classical
$r$-matrix theory for constant $r$-matrices.
\end{abstract}

{\bf  Syst\`emes de Calogero -
 Moser \`a spins associ\'es aux alg\`ebres de Lie simples}
 
{\bf R\'esum\'e}
Nous introduisons des syst\`emes de Calogero - Moser avec spins
 associ\'es aux syst\`emes de racines
des alg\`ebres de Lie simples et en donnons les repr\'esentations de Lax
 (avec param\`etre spectral) associ\'ees ainsi que les crochets de Poisson
fondamentaux. Notre analyse
est bas\'ee sur une approche par la th\'eorie des groupes dont l'esprit est
similaire
\`a la th\'eorie standard
de la $r - $ matrice  classique pour les matrices $r$ constantes.

{\bf Version fran\c caise abr\'eg\'ee}

Le syst\`eme  $\frak{sl}(N)$ rationnel   de Calogero-Moser avec spins 
 a \'et\'e
introduit par Gibbons et Hermsen \cite{GH}
 comme g\'en\'eralisation des syst\`emes
\`a plusieurs particules de Calogero et Moser.
Ainsi qu' en l'absence de spin, il existe aussi des versions trigonom\'etriques
et elliptiques de ces g\'en\'eralisations.
Ces derni\`eres ann\'ees, ces mod\`eles ont fait l'objet d'une certaine
attention
du fait de leur int\'er\^et dans des contextes divers (voir , par exemple,
\cite{P}, \cite{BAB2}, \cite{DP}, \cite{KM}).
Du point de vue des syst\`emes int\'egrables, les crochets de Poisson
fondamentaux entre
les \'el\'ements de matrices de l'op\'erateur de Lax associ\'e $L$ de ces
mod\`eles pr\'esentent
une caract\'eristique inhabituelle \cite{BAB1} \cite{BAB2};
outre les termes usuels contenant la  matrice $ r$ (dynamique) et
l'op\'erateur de Lax $L$,
il apparait un terme d'anomalie dont la pr\'esence est une obstruction \`a
l'int\'egrabilit\'e.
D\`es lors, dans \cite{BAB1} \cite{BAB2},, par exemple, on montre que pour $\frak{sl}(N)$,
 les mod\`eles sont
 "int\'egrables" seulement sous certaines
contraintes additionnelles portant sur les variables de spin.
Il est \`a noter que, 
 bien que le calcul de  matrices $r$ dynamiques ait \'et\'e
effectu\'e pour ces mod\`eles, 
la signification de ces calculs dans la theorie
des groupes reste mysterviense.
\smallskip
Dans cette Note, nous introduisons des syst\`emes de Calogero - Moser avec
spins associ\'es aux syst\`emes de racines des alg\`ebres de Lie simples et
en donnons
les repr\'esentations de Lax (avec param\`etre spectral) ainsi que les
crochets de Poisson fondamentaux.
Notre th\`eme principal est de d\'evelopper une approche dont l'esprit est
similaire \`a
la m\'ethode bien \'etablie de la $r - $ matrice classique pour $r$
constante
\cite{STS}.
 Nous utilisons la classe de $r - $  matrices  dynamiques
avec param\`etre spectral
qui est classifi\'ee dans  \cite{EV}.
 A une telle $ r - $ matrice $r$, nous associons
une application $R$ pour d\'efinir une alg\'ebroide de Lie \`a l'aide du
Corollaire \ref{cor:R}.
Ensuite nous utilisons $R$ \`a nouveau pour d\'efinir une application
de Poisson de l'espace des phases du syst\`eme de spins
dans le dual de cette alg\'ebroide de Lie (voir Proposition \ref{pro:Poisson}).
A l'aide de cette application de Poisson, nous obtenons le Hamiltonien du
syst\`eme de spins,
 sa repr\'esentation
de Lax (avec param\`etre spectral) et les crochets de Poisson fondamentaux
(voir Th\'eor\`eme A, D\'efinition \ref{defi:spin}
 et le reste de la Section 3 pour des
formules explicites).

A ce stade, nous aimerions faire deux remarques importantes.
Premi\`erement, l'application $R$ construite \`a partir de $r$ satisfait
une \'equation dynamique de Yang -Baxter modifi\'ee. (Proposition \ref{pro:Rmdybe} et
Eq. (2)).
Ceci pr\'esente une analogie \'etroite avec le cas non - dynamique.
Deuxi\`emement, dans les cas rationnels et trigonom\'etriques,
les syst\`emes de spins que nous pr\'esentons sont en
correspondance biunivoque avec certains sous - ensembles du
syst\`eme de racines. Nous avons donc autant de syst\`emes de spins que
de tels sous - ensembles.
\smallskip
Il y a plusieurs questions \'evidentes que nous n'abordons pas dans cette Note.
 L'une d'entre elles est l'int\'egrabilit\'e des mod\`eles de spins
au sens de Liouville et la solution pour les flots correspondants au moyen
de notre approche.
Une autre est la r\'eduction aux syst\`emes habituels de Calogero - Moser
associ\'es
aux syst\`emes de racines des alg\`ebres de Lie simples.
Nous consid\'ererons celles - ci et d'autres questions dans des
publications ult\'erieures.

\section{Introduction}
 
     The rational $\frak{sl}(N)$ spin Calogero-Moser system
 was introduced by Gibbons
and Hermsen \cite{GH} as a generalization
 of the many-body systems of Calogero and
Moser.
As in the spinless case, trigonometric and elliptic versions
of this generalization also exist, and in recent years, these models have 
received some attention due to their relevance in a number of areas (see, for
example, \cite{P}, \cite{BAB2}, \cite{DP}, \cite{KM}).
     From the point of view of integrable systems, an unusual feature of these
models is the fundamental Poisson bracket relations between the elements of
the associated Lax operator $L$ \cite{BAB1} \cite{BAB2}. 
 This is because in addition to
the usual terms involving the (dynamical) $r$-matrix and the Lax operator L, 
there is an anomalous term whose presence is an obstruction to integrability.
Thus, in \cite{BAB1} \cite{BAB2},
 for example, it was shown that the $\frak{sl}(N)$ models are
``integrable"
 only when further constraints are imposed on the spin variables.\footnote{Note
that the constrained manifold is not a Poisson manifold in either
\cite{BAB1} or \cite{BAB2}, although the conserved
quantities Poisson commute on it.}  It
should be noted that although dynamical r-matrices
 have been computed for 
these models, however,  the   group-theoretic meaning
of these calculations remains mysterious.
     
In this Note, we introduce spin Calogero-Moser systems associated with
root systems of simple Lie algebras, and give the associated Lax 
representations 
(with spectral parameters) and fundamental Poisson bracket relations.
Our main theme here is to develop an approach which is similar in spirit to
the well-established classical $r$-matrix
 theory for constant $r$-matrices \cite{STS}.
We shall make use of the class of classical dynamical r-matrices with spectral
parameter which is classified in \cite{EV}.  If $r$ is such an $r$-matrix,
 we first use an associated map $R$
 to construct a Lie algebroid.  Then we use it again
to define a Poisson map from the phase space of the spin system into the
dual of the Lie algebroid which we constructed previously.  This Poisson map
is the cornerstone of our theory as it allows us to derive everything about
the corresponding  spin  model.
 Two remarks are in order here.
First, the map $R$ constructed from $r$
 satisfies a modified dynamical Yang-Baxter
equation.  This is in close analogy with the nondynamical case.  Secondly,
note that in the rational and trigonometric case, the spin systems that we
have are  in one-to-one correspondence with some subsets of the root
system.   Thus we have as many spin systems as these special
subsets.
There are several obvious questions which we do not address in this Note.
One of these is the integrability of the spin models in the Liouville sense
and the solution of the corresponding flows.  Another one is the reduction
to the usual Calogero-Moser systems associated with root systems of simple
Lie algebras.  We shall address these and other issues in subsequent 
publications.

\section{General construction}
Let $\frakg$ be a Lie algebra 
and $\frakh \subset \frakg$ an Abelian Lie subalgebra.
In the presentation that follows, we suppose that $\frakg$ is
finite dimensional. In the infinite dimensional case,
it is understood that some suitable modification  has to be made.
  Consider   $T^*\frakh^{*} \times \frakg^*$ as
a vector bundle over $\frakh^{*}$, and define  a bundle map
$a_{*}: T^*\frakh^{*} \times \frakg^*\lon T\frakh^{*}$
by $a_{*}(q, p, \xi )=i^{*} \xi$, $\forall q\in \frakh^{*} ,
p\in \frakh$ and $\xi\in \frakg^*$, where $i: \frakh\lon \frakg$
is the inclusion map.
If  $R$ is   a map from $\frakh^{*}$ to
 $\call (\frakg^* , \frakg )$ (the space of linear maps from
$\frakg^*$ to $\frakg$) which satisfies $R_{q}^* =-R_{q}$ for
each point $q\in \frakh^{*}$ (here, as well as in the sequel,
$R_{q}$ is denoted the linear map in $\call (\frakg^* , \frakg )$
when evaluating $R$ at the point $q$), we
define a bracket on $\gm (T^*\frakh^{*} \times \frakg^* )$ as follows.
For $\xi , \eta \in \frakg^*$  considered as constant
sections, $h\in \frakh $ considered as a constant
one form on $\frakh^{*}$, and $\omega , \theta \in \Omega^{1}(\frakh^{*} )$,
 define
$[\omega, \theta ]=0, \ 
[h , \xi ]=ad_{h}^*\xi , \ 
[\xi ,\eta ]=d(R\xi , \eta )-ad^*_{R\xi }\eta +ad^*_{R\eta }\xi$,
where $ad^*$ denotes the dual of $ad$: $<ad^*_{x}\xi , y>=<\xi, [x, y]>, 
\forall x, y\in  \frakg$ and $\xi \in \frakg^*$.
Then we can  extend  this  to a bracket $[\cdot , \cdot ]$
for all sections in  $\gm (T^*\frakh^{*} \times \frakg^* )$ by
the usual anchor condition.

\begin{pro}
\label{pro:R}
$(T^*\frakh^{*} \times \frakg^* ,  \ [\cdot , \cdot ])$ is a Lie algebroid
  with anchor map $a_{*} $ iff
\begin{enumerate}
\item The operator $R$ is a  function from
 $\frakh^{*} $
to  $\call (\frakg^* , \frakg )^{\frakh}$,
 the space of $\frakh$-equivariant linear map from
$\frakg^*$ to $\frakg$  ($\frakh$ acts on $\frakg$ by
adjoint action and on $\frakg^*$ by coadjoint action).
\item For any $q \in \frakh^{*}$,
the linear map from $\frakg^* \otimes \frakg^* \lon \frakg$ 
defined by
\begin{equation}
\xi \otimes \eta  \lon
[R\xi , R\eta ]+R(ad^*_{R\xi }\eta -ad^*_{R\eta }\xi)
+ \calx_{i^* \xi}(R\eta ) -\calx_{i^* \eta}(R\xi ) +d(R\xi , \eta )
\end{equation}
is independent of $q \in \frakh^{*}$, and is  $\frakg$-equivariant,
where $\frakg$ acts on $\frakg^* \otimes \frakg^*$ by coadjoint action
and on $\frakg$ by adjoint action. Here 
 $\calx_a$ for $a\in \frakh^*$ means the derivative with respect to
$q$  along the constant vector field defined by $a$.
\end{enumerate}
\end{pro}
This is a reformulation of Theorem 4.3.1 in \cite{BK-S} in terms
 of the operator $R$. In fact from Theorem 4.3.1 in \cite{BK-S},
 we know that $(T\frakh^{*} \times \frakg , \ T^*\frakh^{*} \times \frakg^* )$
is indeed a Lie bialgebroid.

In some applications, $\frakg$ admits a non-degenerate ad-invariant pairing
$(\cdot , \cdot )$.  
If $I: \frakg^* \lon \frakg$
is the  induced  isomorphism,  then we have the following

\begin{cor}
\label{cor:R}
The operator $R:  \frakh^{*} \lon \call (\frakg^* , \frakg )$ defines
a Lie algebroid structure on  $T^*\frakh^{*} \times \frakg^*$
if  the condition (1) in Proposition \ref{pro:R} is satisfied,
 and if  $R$ satisfies the modified dynamical Yang-Baxter equation (mDYBE):
\begin{equation}
\label{eq:mdybe}
[R\xi , R\eta ]-  R(I^{-1}[R\xi , I\eta ]+I^{-1}[I\xi , R\eta ] )
+ \calx_{i^* \xi}(R\eta ) -\calx_{i^* \eta}(R\xi ) +d(R\xi , \eta )
=c[I\xi , I\eta ],  \ \ \forall \xi, \eta \in \frakg^* ,
\end{equation}
for some constant $c$.
\end{cor}

Now suppose that $\rho : X\lon T\frakh^* \times  \frakg$ is   a Poisson map,
 where $T\frakh^*\times \frakg $ is equipped with the   Lie Poisson 
structure  induced from the Lie algebroid structure
on $T^*\frakh^* \times  \frakg^*$  constructed  in Proposition 
\ref{pro:R}.   Set
$L=pr_{2}\smalcirc \rho :X\lon \frakg$
and $\tau =pr_{1}\smalcirc \rho :X\lon   T\frakh^* $.

\begin{pro}
\label{pro:A}
 Let $\pi$ denote the Poisson tensor on $X$, and
$\pi^{\#}: T^*X \lon TX$ its induced bundle map. If
$p: T\frakh^*\times \frakg\lon \frakh^*$ is the bundle
projection, and $m=p\smalcirc \rho$, then\\ 
(1). \begin{equation}
\label{eq:LR}
T_{x}L\smalcirc \pi^{\#}_{x} \smalcirc T_{x}^* L = R_{m(x)}\smalcirc ad^*_{L(x)}+
ad_{L(x)}\smalcirc R_{m(x)} +\calx_{\tau (x)}(m(x))R, \ \ \ \forall x\in X.
\end{equation}
(2). If  $f\in C^{\infty}(\frakg )^{G}$,  then under the 
flow $\phi_{t}$ generated by the Hamiltonian $L^* f$ on $X$, we have
the following quasi-Lax type equation:
\begin{equation}
\frac{dL(\phi_{t})}{dt}=[R_{m(\phi_{t}) } (df(L (\phi_{t}))),
 L(\phi_{t})] -\calx_{\tau (\phi_{t})} (R_{m(\phi_{t})}(df(L(\phi_{t}))).
\end{equation} 
\end{pro}

\section{Spin   Calogero-Moser systems}

Let $\frakg$ be a simple Lie algebra and $\frakh\subset \frakg$
be a fixed  Cartan   subalgebra. We shall denote the
Killing form by $(\cdot , \cdot )$ throughout
the section and make the identifications $\frakg^*\cong \frakg$,
$\frakh^*\cong \frakh$.
By definition \cite{EV},   a  classical
dynamical r-matrix with  spectral parameter associated with the
pair $\frakh\subset \frakg$
is  a meromorphic function
$r \, : \, \h^* \,\times \,\complex \, \to \g \otimes \g $
satisfying (a) the zero weight condition; 
(b) the  generalized unitarity condition; (c) the  residue condition:
$ {\mbox{Res}}_{z=0} \, r(q, z)\, = \, \Omega$; and (d) 
the classical dynamical Yang-Baxter equation (CDYBE):
$$
\Alt  (d_\h r) \, + \, [r^{12}(q, z_{1,2}), r^{13}(q,z_{1,3})]\,+ \,
[r^{12}(q,z_{1,2}), r^{23}(q, z_{2,3})]\,
+\,[r^{13}(q, z_{1,3}), r^{23}(q,z_{2,3})] \, = \, 0 \, .
$$
where $z_{i,j}=z_i-z_j$,  and $\Omega \in (S^{2}\frakg )^{\frakg}$
 is the standard Casimir
element.


By $\L\frakg$, we denote the Lie algebra of Laurent series
$\sum_{n=-N}^{\infty}\xi_{n}z^n$ with coefficients
in $\frakg$, which are convergent in some annulus
$0<|z|<r$ (which may depend on the series).
We define $\L\frakg^*$ (the restricted dual) in a
similar fashion.
Associated with each dynamical $r$-matrix $r$,
there is  an operator $R: \frakh^* \lon \call (\L\frakg^* , \L\frakg)$
defined as follows.  Let $r_{0}>0$ be such that
$r(q, z)$ has no  pole in $|z|<r_{0}$
other than $0$. For $0<|z|<\half r_{0}$, we set
\begin{equation}
\label{eq:rR}
(R_{q} \xi)(z)=p.v. \frac{1}{2\pi i} \oint_{C} (r(q, w-z), \xi (w)\otimes 1)dw
\end{equation}
where $C$ is the circle of radius $|z|$
 with positive orientation,
and $p.v.$ denotes the principal value of the  improper integral. Note that
$r(q, w-z)$, as a function of $w$, has no singularity in the interior
of $C$, and the principal value of the integral
can be shown to exist. 

By some calculation using analysis, we prove
the   following:

\begin{pro}
\label{pro:Rmdybe}
 The operator $R$ defined by Equation (\ref{eq:rR})
is  in $\call (L\frakg^* , L\frakg)^{H}$ and 
satisfies the mDYBE (Equation (\ref{eq:mdybe})) with $c=-\frac{1}{4}$.
\end{pro}

According to Corollary \ref{cor:R},  we can use $R$ to equip
  $T^*\frakh^{*} \times L\frakg^*$
 with a Lie algebroid structure, and therefore 
$T\frakh^{*} \times L\frakg$ admits  the Lie Poisson structure.
On the other hand, consider   $T^*\frakh^{*}$ with 
the  canonical cotangent  symplectic structure, $\frakg^*$
     with the  plus Lie Poisson structure, and equip $T^*\frakh^{*} \times \frakg^*$
with the
product Poisson structure.
Then we have 

\begin{pro}
\label{pro:Poisson}
The map $\rho : T^*\frakh^{*} \times \frakg^*\lon T\frakh^{*} \times L\frakg$ 
given by
$$(q, p, \xi )\lon (q, -i^*\xi, p+r_{-}^{\#}(q)\xi), \ \ q\in \frakh^{*},
p \in \frakh, \xi \in \frakg^* , $$
is a Poisson map, where $i: \frakh \lon \frakg$ is the natural
inclusion, and $r_{-}^{\#}(q): \frakg^*\lon L\frakg$ is
defined by $((r_{-}^{\#}(q)\xi)(z),  \eta )=(r(q, z), \eta \otimes \xi)$
for all $\xi, \eta \in \frakg^*$.
\end{pro}

Combining Proposition \ref{pro:Poisson} with
Proposition \ref{pro:A}, we are lead to the main theorem of this Note. 

{\bf Theorem A}
{\em Assume that $r$ is a     classical
dynamical $r$-matrix with  spectral parameter. Then
$L:  T^*\frakh^{*} \times \frakg^*\lon  L\frakg$, 
$(q, p, \xi )\lon  p+r_{-}^{\#}(q)\xi $
satisfies the  fundamental Poisson bracket relations:
\begin{eqnarray}
\{L_{1} (z), L_{2}(w)\}&=&-[r^{12}(q, z-w), L_{1}(z)]+[r^{21}(q, w-z), L_{2}(w)]
-\calx_{i^*\xi}r(q, z-w)\\
&=&-[r^{12}(q, z-w),  L_{1}(z)+ L_{2}(w) ]-\calx_{i^*\xi}r(q, z-w)
\end{eqnarray}
Moreover, if $f\in C^{\infty}(L\frakg)$ is ad-invariant, then under  the  
flow  $\phi_{t}$ generated by the Hamiltonian $L^* f$, we have the
following   quasi-Lax type equation:
\begin{equation}
\frac{dL(\phi_{t})}{dt}=
[R_{m(\phi_{t})} (df(L (\phi_{t}))), L(\phi_{t})] 
+\calx_{i^{*}(pr_{2} \phi_{t} )} (m(x)) (R_{m(\phi_{t})}(df(L(\phi_{t}))).
\end{equation}
}

Let $E$ be the quadratic function
\begin{equation}
E(\xi)=\half \oint_C (\xi(z), \xi (z))\frac{dz}{2\pi i}, \ \ \ \forall \xi\in L\frakg,
\end{equation}
where $C$ is a small circle around the origin.
Clearly, $E$ is an ad-invariant function on $L\frakg$.

\begin{defi}
\label{defi:spin}
Assume that $r$ is a     classical
dynamical $r$-matrix with  spectral parameter.
The Hamiltonian system on $ T^*\frakh^{*} \times \frakg^*$  generated
by the Hamiltonian function:
\begin{equation}
H(q, p, \xi)= (L^*E)  (q, p, \xi)= \half \oint_C (L(q, p, \xi), L(q, p, \xi))
\frac{dz}{2\pi i}
\end{equation}
is called the  spin Calogero-Moser system associated to the
dynamical r-matrix $r$.
\end{defi}

In  \cite{EV}, Etingof and Varchenko 
obtained a complete classification of classical
dynamical $r$-matrices. Up to  gauge transformations, they obtained
canonical forms of the three types  (rational, trigonometric
and elliptic) of dynamical $r$-matrices.
 For  each of these dynamical r-matrices,
one can associate a spin Calogero-Moser system on $ T^*\frakh^{*} \times \frakg^*$.  We will  list all of them below.
First, let us fix some notations. Let $\frakg =\frakh \oplus 
\sum_{\al \in \Delta}\frakg_{\al}$ be the root
space decomposition. For any positive
root $\alpha \in \Delta_{+}$, fix basis  $e_{\al}\in \frakg_{\al}$
and $e_{-\al}\in \frakg_{-\al}$  which are dual with respect to
$(\cdot , \cdot )$.
Fix also  an orthonormal  basis $\{h_{1}, \cdots , h_{N}\} $
of $\frakh$, and write $p=\sum_{i=1}^N p_{i}h_{i}$ and $\xa= (\xi , e_{-\al})$,
 for $p\in \frakh$ and $\xi \in \frakg^*$.

{\bf I. Rational case}

\be
&&r(q, z) \, = \,  {\Omega  \over z} \,
+ \sum _{\alpha  \in \Delta'}  \,
{1\over (\al, q) } \,
e_\al \otimes  e_{-\al};\\
&&H(q, p, \xi )=\half \sum_{i=1}^{N}p_{i}^{2}-\half \sum_{\alpha \in \Delta'}
\frac{\xa \ax}{(\alpha , q)^2},\\
&&L(q, p, \xi )(z)=p+\frac{\xi}{z}+\sum_{\alpha \in \Delta'}
\frac{\xa}{(\alpha , q)}e_\al ,
\ee
where  $\Delta' \subset \Delta$ is a set of roots
 closed with respect to the addition and multiplication by $-1$.

{\bf II.  Trigonometric  case}

\be
&&r(q, z)\,=\, (\mbox{cot}\, z +\third z ) \,\sum_{i=1}^N h_i\otimes h_i\,
+ \, \sum _{\al \in \Delta (\Pi' )} \,{ \mbox{sin}\, ((\al, q) + z)
\over
 \mbox{sin}\, (\al, q) \,  \mbox{sin}\, z} 
e^{\third z (\alpha , q)}\,e_\al \otimes  e_{-\al} \\
&&\ \ \ \ \ \ \ \ \ \ + \, \sum _{\al \in \Delta_{+}-\Delta (\Pi'  )}\,{ e^{-iz}
\over
 \mbox{sin} \,z} e^{\third z (\alpha , q)} \,e_\al \otimes  e_{-\al} 
+ \, \sum _{\al \in \Delta_{-}-\Delta (\Pi'  )}\,{ e^{iz}
\over
 \mbox{sin}\, z}e^{\third z (\alpha , q)}\,e_\al \otimes  e_{-\al}, \\
&&H(q, p, \xi)=  \half \sum_{i=1}^{N}p_{i}^{2}-\half
\sum _{\al \in \Delta (\Pi'  )}(\frac{1}{\sin^{2}(\al, q)}-\frac{1}{3})\xa \ax
-\frac{5}{6}\sum _{\al \in \Delta-\Delta (\Pi'  )}\xa \ax\\
&&L(q, p, \xi)(z)=p+(\cot{z}+\frac{1}{3} z)(i^* \xi)+
\sum _{\al \in \Delta (\Pi'  )}{ \mbox{sin}\, ((\al, q) + z)
\over
 \mbox{sin}\, (\al, q) \,  \mbox{sin}\, z}e^{\frac{1}{3}z (\al , q)}\xa e_\al \\
&&\ \ \ \  \ \ \ \ \ \  \ \ + \sum _{\al \in \Delta_{+}-\Delta (\Pi'  )}\,{ e^{-iz}
\over
 \sin{ \,z}} e^{\frac{1}{3}z (\al , q)}\xa e_\al
+\sum _{\al \in \Delta_{-}-\Delta (\Pi'  )}\,{ e^{iz}
\over
 \mbox{sin}\, z}e^{\frac{1}{3}z (\al , q)}\xa e_\al .
\ee
Here $\Delta =\Delta_{+}\cup \Delta_{-}$
is a polarization of $\Delta$, 
 $\Pi' $ is a  subset of the set of simple  roots,
and  $\Delta (\Pi'  )$
denotes the set of all roots which are linear
combinations of roots from $\Pi' $.

{\bf  III.  Elliptic case }

\be
&&r(q, z,)\,=\, \zeta (z)\sum_{i=1}^N h_i\otimes h_i\,
- \,\sum_{\al \in \Delta} l(q, z ) e_\al \otimes e_{-\al} \\
&&H(q, p, \xi)=  \half \sum_{i=1}^{N}p_{i}^{2}-\half
\sum _{\al \in \Delta }\calp ((\al , q ))\xa \ax\\
&&L(q, p, \xi)(z)=p+\zeta (z)(i^*\xi ) -\sum _{\al \in \Delta } l((\al , q ), z)
\xa e_{\al}
\ee
where $\zeta (z)=\frac{\sigma' (z)}{\sigma (z)}$, $\calp (z)
=-\zeta' (z)$,
 $l(q, z)=-\frac{\sigma (q+z)}{\sigma (q)\sigma (z)}$, and
$\sigma (z)$ is the Weierstrass $\sigma$ function of periods 
$2\omega_{1}, 2\omega_2$. 

\noindent
{\bf Theorem B}
{\em  The Hamiltonians of the spin Calogero-Moser
systems are invariant under the canonical action
$H\times  T^*\frakh^{*} \times \frakg^*\lon T^*\frakh^{*} \times \frakg^*$,
$x\cdot (q, p, \xi )=(q, p, Ad_{x^{-1}}^*\xi )$ with 
momentum map $J: T^*\frakh^{*} \times \frakg^* \lon \frakh^{*}$
given by $J(q, p, \xi ) =i^*\xi$. If $\Sigma$ denotes
the set defined by $\calx_{i^*\xi}R=0$, then
$\Sigma=J^{-1}(0)$ in the trigonometric and elliptic
cases, while $\Sigma =J^{-1} ((\Delta ')^{\perp})$ in
the rational case. Thus in each case, $\Sigma$ is invariant
under the dynamics and we have $\frac{dL}{dt}=[R (M) , L]$
on $\Sigma$, where $M(q, p, \xi )(z) =L (q, p, \xi )(z)/z$. }

{\bf Acknowledgments.}  In addition to the funding sources mentioned
in the first footnote, we like to thank several institutions
for their hospitality while work on this project was being done:
MSRI (Li),  and Max-Planck Institut (Li and Xu). We  also thank
Serge Parmentier 
 for his help in preparing the French abbreviated version.

\end{document}